\documentclass[12pt,reqno]{amsart}

\usepackage{amssymb}
\usepackage{ifthen}
\usepackage{verbatim}
\usepackage{epic,eepic}

\renewcommand{\phi}{\varphi}

\newcommand{\Z}{{\mathbb Z}}

\newcounter{step}

\theoremstyle{plain}

\theoremstyle{definition}

\begin{document}

\title{Structure Theory of Set Addition III. Results and Problems.}

\author{Gregory A. Freiman}
\address{G. A. Freiman\vspace{-0.25cm}}
\address{School of Mathematical Sciences, Tel Aviv University,
Tel Aviv 69978, Israel.}
\email{~grisha@post.tau.ac.il}

\date{}



\maketitle

\section{Introduction}

\noindent
This review is motivated by a powerful and breathtaking
development in ``Additive Combinatorics", the direction of
study connected with the names of
Timothy Gowers, Jean Bourgain, Terence Tao and Ben Green.

The results and methods of Inverse Additive Number Theory
have been used substantially in this new field.

Additive Combinatorics now looks towards algebra and
computation.  Results in these areas existed in Inverse
Additive Number Theory and they will be discussed below.

This paper is a continuation of the reviews [32] and  [33].

\section{The main theorem}

\noindent
We consider a finite subset $A$ of $\Z^n$ (and often $\Z$)
with cardinality $|A|=k$. 
We define 
$$2A=A+A=\{x: x=a+b, a\in A, b\in A\}.$$
We say that $A$ has the small doubling property if
$$|2A|<Ck,$$
where $C$ is a positive constant or a slowly increasing
function.
The number $$\frac{|2A|}{k}$$
is called the doubling coefficient of $A$. 

We define the notion of additive isomorphism of
$A \subseteq \Z^n$
onto $B \subseteq \Z^m$ as a bijection $\phi:A\rightarrow B$
such that, for all $a,b,c,d$ in $A$, we have 
\begin{align}
\phi (a)+ \phi (b)=\phi (c) +\phi (d)  
\end{align}
if and only if 
$$a+b=c+d.$$

We write $A\sim B$ and note that $\sim$ is an equivalence
relation.
We name a \hbox{$d$-dimensional} parallelepiped the subset
$D\subseteq \Z^d$
$$D=\{(x_1,x_2,\ldots,x_d: x_i\in \Z,\,0 \le x_i <
h_i,\,h_i \ge 2,\, 1 \le i \le d)\},$$
where
$$|D| =h_1 h_2\dots h_d.$$

\medskip
\noindent
{\bf Main Theorem.}
{\it For all finite subsets $A \in \Z$ for which $|2A|<Ck$,
where $C$ is a positive constant or $C=C(k)$ is a slowly
increasing function, there exist $c=c(C)$ and a
parallelepiped $D$ with $d \le [C-1]$ and $|D| <ck$ such
that $A \subseteq \phi (D)$, where $\phi$
is an isomorphism.} \hfill $\square$

\medskip
For a given $A\in\Z$, we consider the images $\phi (A)$
under all isomorphisms $\phi: A \rightarrow \Z^n$, for all
$n \ge 1$. The dimension of any $\phi (A)$ 
is the dimension of the smallest affine subspace
containing it.
The maximum dimension of $\phi (A)$ for all such $\phi$ is
called the {\it dimension} of $A$ and denoted $d(A)$. 

We will define the notion of {\it volume} of $A$.
Let $A \in \Z^n$, where $n$ is the dimension of $A$.
For each isomorphism $\phi: A \rightarrow \Z^n$,
let us take the convex hull of $\phi(A)$ and let  $|\phi(A)|$
be the number of integer points in it. Then we define
$$V(A)=\min_{\phi}|\phi(A)|.$$

The proof of the Main Theorem was obtained gradually 
(see [14]--[19]). 
The final text of the proof may be found in [13] and
Yuri Bilu's paper [2]. 
Somewhat different proofs were given by Ruzsa [41], [42],
Mei-Chu Chang [11] and Sanders [43], with better estimates
of volume. 

Now, some remarks on estimates of $V(A)$, the volume of $A$.
We can also estimate $V(A)$ by the volume
$$|D|=h_1h_2\dots h_d$$
of a $d$-dimensional parallelepiped $D$ with a minimal $|D|$
and $A \subseteq \phi(D)$, where $d$ is the dimension of $A$
and $\phi$ is an isomorphism.
The set $\phi(D)\subset \Z$ is called a generalized 
proper arithmetic progression of dimension $d$. 

Existing estimates were obtained with the help of $\phi(D)$,
where $\phi$ is not an isomorphism but a homomorphism giving
a generalized arithmetic progression which is not proper. 
As a result, we obtain an estimate of the volume $V(A)$
and dimension of $A$, which wait to be improved.
Below, we will give a hypothetical best possible estimate.

\section{Hypothetical value of $V(A)$}

\noindent
We begin with an example which does not cover all possible
values of $T=|2A|$ but is very simple and representative. 

Let
 $$A=\{0,1,\ldots,m-1, 2(m-1), 4(m-1),\ldots,2^{c-2}(m-1)\}$$
where, as usual, $|A|=k$,
and
\begin{align}\label{e:3}
3\le m \le k,  
 \end{align}
 \begin{align} 
 c=k+2-m,
 \end{align}
so that $2 \le c \le k-1$.

   In this case we have $d(A)=1$ and the values of $V$ and $T$
have a very simple form for $c$ which is not very large:
 \begin{align}\label{e:V}
  V=2^{c-2}(m-1)+1=\frac{1}{4}2^c(k+1-c)+1=\frac{1}{4}2^ck+O(c2^c)
 \end{align}   
and 
\begin{align}\label{e:def-T}
T&=(k-1)+(k-2)+\dots+m+2m-1=\frac{k-1+m}{2}(k-m)+2m-1=
\nonumber\\
&=\frac{k-1+k+2-c}{2}(c-2)+2(k+2-c)-1=\nonumber\\
&=ck-\frac{(c-1)(c-2)}{2}-2(c-2)-1=ck-\frac{c^2+c-4}{2}=\nonumber\\
&=ck+O(c^2).
\end{align}
If $c=O(\log^\delta k)$, $\delta <1$, then, for $V$ and $T$,
(4) and (5) give very simple expressions, sufficient for
today's use.

Now, we will give examples of the sets $A$, for which $T$
takes all possible values. We will see, that the ``net" of
sets $A$, given earlier, represents the general situation
in a very good way.

For values $m$ and $c$ as in (2) and (3), let
\begin{align}\label{e:a_m}
a_{m-1}=m-1+b,
\end{align}
where
\begin{align}\
0\le b \le m-3=k-c-1,
\end{align}
and denote
\begin{align}\label{e:A}
A=\{0,1,\ldots,m-2,a_{m-1}, 2a_{m-1},
4a_{m-1},\ldots,2^{c-2}a_{m-1}\},
\end{align}
where conditions (6) and (7) are valid. We get
\begin{align}\label{e:V}
V&=2^{c-2}a_m+1=\frac{1}{4}2^c(m-1+b)+1=\nonumber\\
&=\frac{1}{4}2^c(k-c+1+b)+1=\frac{1}{4}2^c(k+b)+O(c2^c),
\end{align}
and 
\begin{align}\label{e:T}
T=(k-1)+(k-2)+\dots+m+2m-1+b=ck+b+O(c^2).
\end{align}

If $d(A)>1$, we build the set $A$ in the following manner.
Let
\begin{align}\label{e:A}
A=A_d=\{A'e_1,e_2,\ldots,e_d\}
\end{align}
where $A'\subseteq \Z$ will be defined a bit later.

For a given value $T$, we define
\begin{align}
  |2A'|=T-\big(k+(k-1)+(k-2)+\dots+(k-d+2)\big).
 \end{align} 
We have 
\begin{align}\label{e:k'}
k'=|A'|=k-d+1.
\end{align}
For values $k'$ and $|2A'|$ given in (12) and (13), 
we can build the set $A'$ as in (8)
and the set $A_d=A$ with the help of (11). 

In the case where the values of $k$ and $T$ are
given,
we can also get the corresponding values of $V(A)$ as a 
function of $k$ and $T$. It was done in [13, p.\,37],
where the hypothesis about the best estimate of $V(A)$ was
formulated. 
We may reproduce it now in a nicer form.

\medskip
\noindent
{\bf Hypothesis.} {\it The inequality
\begin{align}\label{e:VA}
V(A) \le \frac{1}{4}2^c(k-c+1+b)+1
\end{align}
gives an estimate from above for $V(A)$ 
for a given $k$ and $T$ which cannot be improved.}

\medskip
If this hypothesis were proved to be true, the proof
would be important not only as an improvement of an
estimate of $V(A)$ for a given doubling coefficient $T$,
but also the range of values $T$ would cover all
possible values at the time when existing estimates make
sense only for rather small values of $T$.
     
The set $A$ may, because of (11), be put in a
parallelepiped $D$ with edges 
$$h_1=V(A'),\quad h_2=h_3=\dots h_d=2$$
and 
$$V(D)=2^{d-1}V(A').$$
We may define the volume of a set $A$ with the help
of the
volume of $D\supseteq A$ for a minimal $|D|$ for such $D$ 
(a proper $d$-dimensional arithmetical progression) and 
formulate a hypothesis about minimal volume 
in a simpler though less exact way.

\section         {On a family of extremal sets}
        
\noindent
We would  now like to go one step deeper into the study    
of the structure of $A$. We have used a concrete example 
of $A$ with given values $k$ and $T$
for which the maximal value given in (14) is achieved.
   
Let us ask about a list of such $A$. 
Here, as is usual in solving an inverse additive
problem, see [32, p.\,5],
we study the structure of the set when its
characteristic is equal to its extremal value.

As an example of such sets, let us build a family of
sets $\{B\}$ which are ``friendly"
to $A$ defined in (8) in the following manner:
\begin{align*}
&|B|=|A|=k,\\
&B=\{b_0,b_1,\ldots,b_{k-1}\}.
\end{align*}
Define $B(s)$ as a subset of $B$ containing the first $s$
elements of $B$:
$$B(s)=\{b_0,b_1,\ldots,b_{s-1}\}.$$
Let $B(m)$ be one of the sets described in Theorem 1.9 of
[13], and in the formulation given at the end of p.11,
we have 
$$b_0=0,b_{m-1}=a_{m-1}=m-1+b,\quad |2B(m)|=2m-1+b.$$
Let $p_m$ be the number of such sets. 
Define the sets $B(s)$ for $m\le s\le k$ in the 
following manner: Let the sets for some value which is 
less than or equal to some chosen value of $s$ be already
known.
Take $B(s)$ which is one of them.
Define 
$$\operatorname{sym}(B(s))=\max B(s) -B(s),$$
and, for a chosen $B(s)$, take two possibilities 
\begin{align}
B'(s+1)=\big\{B(s)\cup\{2^{s-m+1}a_{m-1}\}\big\}
\end{align}
or 
\begin{align}
B''(s+1)=\{\operatorname{sym}B'(s+1)\}.
\end{align}
Note that for $s>m$ we have 
$$\max B(s)=2^{s-m}a_{m-1}$$
and so
\begin{align}
\max B(k)=2^{k-m}a_{m-1}=1/4*2^c(k-c+1+b).
\end{align}
The number of sets in the family $\{B\}$ is equal
to $2^{c-2}p_m$.

     The values $V$ and $T$, for each set $B$  
just built, are the same as in (9) and (10).

Let us give a numerical example.
Take $k=8$, $m=6$ and because of (3) we get $c=4$.
Further, take in $a_{m-1}=m-1+b$, $b=2$ and $a_5 =7$.

For sets $\{B(6)\}$ we obtain the following list:
\begin{multline}
\{B(6)\}=\{0,1,2,3,4,7\},\{0,1,2,3,5,7\},
\{0,3,4,5,6,7\},\\
\{0,2,4,5,6,7\},\{0,2,3,4,5,7\}.
\end{multline}
Therefore, the number of these sets is $p_6=5$.

From (9) and (10) we get
\begin{align}
T=7+6+12-1+2=26,
\end{align}
\begin{align}
V=\frac{1}{4}2^4(6-1+2)+1=29.
\end{align}

Take, for example, the set $\{0,1,2,3,4,7\}$ from
$\{B(6)\}$.
Using (9) and (10), we get two sets $\{0,1,2,3,4,7,14\}$ and
 $\{0,7,10,11,12,13,14\}$ from $\{B(7)\}$
  and four sets $\{0,1,2,3,4,7,14,28\}$,
  $\{0,14,21,25,26,27,28\}$, 
  $\{0,7,10,11,12,13,14,28\}$
   and\break
  $\{0,14,15,16,17,18,21,28\}$ from $\{B(8)\}.$
  
  As a result we will get 20 sets from $\{B(8)\}$.
  We gave a partial list of extremal sets which gives some 
  idea of how to construct them.

\section {Partial case: $T \le 4k-7$}

\noindent
Let us use the following lemma, see [13, p.\,24].

\medskip
\noindent
{\bf Lemma.} {\it Let $A \subseteq E^n$ be a finite set of
dimension n with $|A|=k$.
 Then $$T\ge(n+1)k-\frac{n(n+1)}{2}.$$
}

\medskip
We see that, if $n=2$, then 
\begin{align}
T \ge 3k-3
\end {align}
 and if $n=3$ then
 \begin{align}
T \ge 4k-6.
  \end{align}
From (21) and (22), we see that if $T \le 3k-4$ then $d(A)=1$.
This case was studied in detail in [13, p.\,11--14].
As to the case,                             
$$3k-3 \le T \le 4k-7,$$
it gives a good example of the Hypothesis for a very simple
case.

    In this case, (8) takes the form
   $$ A={0,1,2,\ldots,k-3,k-2+b,2(k-2+b)}$$
   where
   $$T=k-1+2(k-1)-1+b=3k-4+b, $$
$$  0\leq b \leq k-3.$$ 

Then, according to the hypothesis, $A$ is part of a segment
of length $2k-3+2b$ or has dimension two.
   This case is much simpler than the general case
(on the second page of this article); I had formulated it
immediately after Theorem 1 was proved, but even in this
 simpler case, the problem is still open (for some
progress, see Jin [38]).
    
\section    {Group theory}
    
\noindent
Over the last few years, there has been exciting development
in group theory (Tao [44], Breuillard and Green [7], [8], 
Green [36], Fisher, Katz and Peng [12], Hrushovski [37])   
    connected with sets with small doubling.
I'm referring to the results in approximation groups.

An approximation group is a finite subset of a group which,
together
with the property of small doubling, contains a unit and
has a symmetry property.

Here are some of my old results in this direction:

\medskip
\noindent{\bf Theorem I.} 
{\em Let $G$ be a torsion-free group, $K$ and $M$ subsets
of $G$, for which
$$2\leq|K|,|M| $$
 and
     $$|KM|=|K|+|M|-1, $$
     then $K$ and $M$ have the form 
 \begin {align}
     K=\{a,aq,\ldots,aq^{l-1}\} 
 \end {align}                        
and 
     $$M=\{b,qb,\ldots,q^{h-1}b\},$$
 where
      $$a,b,q \in G,q\neq 1.$$}
      
\medskip
      And what about $K=M$?

\medskip
\noindent{\bf Theorem II}.
{\em If $|K^2|=2|K|-1$ then $K$ has the form $(23)$ and either
$aq=qa$ or $qaq=a$ holds.}
      
\medskip
      The proof is in [28].
      
The immediate problem here is to describe the structure of
the subset $A\subset G$ in a case when $|A^2|\leq 3k-4$
and the group $G$ is torsion-free.
In the study of approximate groups, the study of subsets of
groups of special kinds take place, 
and this is a very natural and rewarding
turn.  
If $G$ is not torsion free then the doubling coefficient may
be less than two.
In this case, the structure of $A$ was described in [20] for
the doubling coefficient 8/5. What about bigger coefficients?
     
\section{All small subsets have small doubling}
     
\noindent
Now, we will describe the study of group structure, which is
different from the fast growing study of approximate
groups, but which also uses the idea of small doubling
as a main tool, and sometimes arrives at similar results. 
They show that some subsets of a group have the structure
of a subgroup, a coset or a union of cosets for some normal
subgroup.
Let us translate the notion of isomorphic subsets to
the notion of similarity of multiplication tables
(Latin squares).
We will now call an algebraic operation a multiplication, and
take a product $ab$ instead of the sum $a+b$.
We will enter he elements of the sets $D$ and $F$ 
into the following tables:
\begin{align}     
\begin{tabular}{l l l l}
  ~ & a & b & c \\
  a & A & B & C \\
  b  \\
  c  \\
\end{tabular} 
\end{align} 
and
\begin{align}     
\begin{tabular}{ l l l l}
  ~ & x & y & z \\
  x  \\
  y  \\
  z  \\
\end{tabular} 
\end{align} 
Here $a,b,c,\ldots$ in $D$ and $x,y,z,\ldots$ in $F$.  

Entries in a multiplication table may be put in some
freely chosen order.
So, let us change the order of columns in (25), and 
correspondingly the order of lines in (25), in such a way
that we will get the following table instead
of (25):
\begin{align}     
\begin{tabular}{ l l l l}
  ~ & $\phi (a)$ & $\phi (b)$ & $\phi (c)$ \\
  $\phi (a)$  \\
  $\phi (b)$  \\
  $\phi (c)$  \\
\end{tabular} 
\end{align}  

As to notation inside the multiplication table (24), the
elements of $D^2$ and $D$ may be different, and we
denote elements of $D$ by lower-case Latin letters and
elements of $D^2$ by upper-case ones, and we denote equal
elements of $D^2 $ by the same letter and different ones
by different letters.
And now, how is condition (1) in multiplicative form
     $$ab=cd     \Leftrightarrow  \varphi(a) \varphi (b)= \varphi(c) \varphi(d) $$
formulated in the case of Latin squares?
The condition $ab=cd$ means the element of the table on the
line with entry $a$ and on column $b$ is equal to the
element with line $c$ and column $d$. 
The condition $\varphi(a)\varphi(b)=\varphi(c)\varphi(d)$
tells us that elements in the same places in the table (26)
are equal and have the same notation. If the equality is not
valid in table (25), it is not valid for corresponding
places in table (26), as well.
So, as formulated in [22, p.\,141], let $D$ and $F$ be
Latin squares of the same order; also, let $\theta$ be a
bijection of $D$ onto $F$.
Applying $ \theta $ to $D$, we have a new Latin square
$\theta[D]$.
If $F$ can be obtained from $\theta[D]$ by a permutation of
rows and the same permutation of columns of
$\theta[D]$, then $D$ and $F$ are said to be isomorphic.
In paper [22] the classification of all isomorphic classes 
  of three-element sets was presented.
The number of commutative classes is six, and of
noncommutative classes is 45.              
Having these and similar results in mind, we can now ask 
plenty questions.

Some examples.

1) Let us have, for each three-element sets of a
group $G$, $|A^2|\leq a$, $a\leq 8$.            
    All such groups were described in [1], [39].
    What if $|A|=4$ and $|A^2|\leq a$, $a\leq 15? $
    
\smallskip
2) Take one subset of three elements with a given
multiplication table. Describe all groups without
this Latin square.
Do the same for a given family of Latin squares.
Looking at groups in an artistic way, we can think of a
group as a complete building. The subgroups will then
be, say, storeys of this building, and small subsets,
described with the help of their multiplication
tables, as the bricks used to construct the building. 
The kind of material used in the course of construction
has an influence on the properties of the whole building.
We may begin to analyze kinds and numbers of different kinds
of small subsets which are used in groups of some known
kinds, say, crystallographic groups.
It would be important,  though not that simple, to find some
connection between the properties of such groups and some 
properties of their composition.
    
 \section   { Subsets of special elements}
    
\noindent
Let us introduce a notion of an $(n,m)$-special element in
a group $G$.
Element $a\in G$ is called $(n,m)$ special if
$|K^n|\leq m< 2^n$, for all  $K=\{a,g\}$, $g\in G$,
where $n,m$ are positive integers $m,n\geq 2$.
By $K^n$ we mean the set of elements $g\in G$ having at least
one representation of the form
$g=k_1,k_2,\ldots,k_n$.
The set of all the $(n,m)$-special elements is denoted by
$S_{m,n}(G)$.
The notion of a special element is attributed  to J.G.\ 
Berkovich.
In [6], it was shown that $S_{2,3}(G)$ and
$S_{3,5}(G)$ are
normal subgroups of $G$.
Some further results and questions in this direction
can be found in [6].
The main problem here: Let $n$ be fixed. For which $m$ the
set $S$ is a normal subgroup of $G$,union of cosets,...?
   
 \section  {Computation complexity}
   
\noindent
Results in Integer Programming were obtained mainly
between the years 1987 and 1990, see [9], [10],
   [21], [25], [26], [27], [30], [31].
    The main problem which was studied there:
    Let an equation
    $$     a_1 x_1+a_2 x_2+\dots +a_m x_m =b     $$
be given, where $a_i$ are different large positive integers.
Do solutions exist of this equation or, perhaps, there are no
solutions?
    The unknowns $x_i$ take the values 0 \hbox{or 1}.
Direct computation, see, for example [40], gave a solution to
the problem for $m$ equal to two-three hundred, and $a_i$ of
order $10^{10}$.
The use of analytical methods of Additive Number Theory
enabled us to formulate general conditions when a solution
and, in fact, many solutions exist.
Methods of Inverse Additive Number Theory enabled us to build
the structure of the set $A={a_i}$ in the case where there
are no solutions. 
Algorithms and programs were built for values of $m$ of 
order $10^7$ allowing us to produce computations in a
split second (see [10]). 
A problem for further study is to go from one equation to
a system of several equations (see [31] and the thesis of
Alain Plagne).
 


\bigskip
\noindent
{\bf References}
\medskip
\begin{enumerate}

\item    Ya.G. Berkovich, G.A. Freiman, C. Praeger \\
       Small squaring and cubing properties for finite groups. \\
       Bull. Australian Math. Soc. 44, n.3, 429-450 (1991).\\
              
\item Yuri Bilu\\
      Structure of sets with small sumset.\\
      Asterisque 258,1999,77-108.  \\ 

\item    L.V. Brailovsky, G.A. Freiman\\
       On a product of finite subsets in a torsion-free group.\\
       J. of Algebra 130, 462-476 (1990).\\
   
\item  L.V. Brailovsky, G.A. Freiman \\
       On two-element subsets in groups. \\
       Ann. of the N. Y. Acad. Sci. 373, 183-190 (1981).\\

\item  L.V. Brailovsky, G.A. Freiman \\
       Groups with small cardinality of the cubes of their two
       element subsets.\\
       Ann. of the N.Y. Acad. Sci. 410, 75-82 (1983).\\
           
\item    L. Brailovsky, G. Freiman, M. Herzog\\
       Special elements in groups.\\
       Proceedings of the Second International Group Theory
       Conference, Supplemento ai Rendiconti del Circolo Matematico di
       Palermo, Bressanone-Brexen, June 11-17, 1989, Ser II, no. 23,
       33-42 (1900).\\
       
\item Emmanuel Breuillard, Ben Green  \\
      Approximate groups,I:the torsion-free nilpotent case.\\
      arXiv:0906.3598v1 [math.CO],19 Jun 2009.\\

\item Emmanuel Breuillard, Ben Green\\
      Approximate groups, II: the solvable linear case.\\
      arXiv:0907.0927v1 [math.GR],6 Jul 2009.\\
   
\item    M. Chaimovich, G. Freiman, Z. Galil \\
       Solving dense subset-sum problems by using analytical number
       theory.\\
       Journal of Complexity 5, 271-282 (1989).\\
       
\item   M.Chaimovich.
      Subset-Sum Problems with different summands:Computation,
      Discrete Applied Mathematics,27,(1990),277-282.\\

\item M.-C.\ Chang\\
      A polynomial bound in Freiman's theorem.\\
      Duke Math.J.113(3),2002, 399-419.\\

\item D.Fisher,N.H.Katz and I.Peng\\
      On Freiman's theorem in nilpotent groups,\\
      http://arxiv.org./abs/0901.1409.\\

\item G.A. Freiman \\
      Foundations of a Structural Theory of Set Addition. \\
      Translations of Mathematical Monographs, Vol. 37. Amer. Math. Soc.
       Providence,.I.R.,108 pp. (1973).\\
      
\item  G.A. Freiman\\
       The addition of finite sets. I\\
       Izv. Vyss. Ucebn. Zaved. Matematika 6(13), (Russian) 202-213 (1959).\\

\item    G.A. Freiman\\
       Inverse problems in additive number theory.  On the addition
of sets of  residues with respect to prime modulus.\\
       Dokl. Akad. Nauk SSSR 141, (Russian) 571-573 (1961).\\

\item    G.A. Freiman\\
       Inverse problems of additive number theory VII.  Additions of
       finite sets IV.  The method of trigonometric sums.\\
       Izv. Vyss. Ucebn. Zaved. Matematika 6(31), (Russian) 131-144
       (1962).\\

\item    G.A. Freiman \\
       On the addition of finite sets.\\
       Dokl. Akad. Nauk SSSR 158, (Russian) 1038-1041 (1964).\\

\item    G.A. Freiman\\
       Inverse problems of additive number theory VIII.  On a
       conjecture of P. Erdos.\\
       Izv. Vyss. Ucebn. Zaved. Matematika 3(40), (Russian) 156-169
       (1964).\\

\item    G.A. Freiman\\
       Inverse problems of additive number theory IX.  The addition of
       finite sets V.\\
       Izv. Vyss. Ucebn. Zaved. Matematika 6(43), (Russian) 168-178
       (1964).\\

\item    G.A. Freiman \\
       Groups and the inverse problems of additive number theory. \\
       Kalinin. Gos. Univ. Moscow, (Russian) 175-183 (1973).\\

\item    G.A. Freiman \\
       An analytical method of analysis of linear Boolean equations. \\
       Ann. of the N. Y. Acad. Sci. 337, 97-102 (1980).\\

\item    G.A. Freiman \\
       On two and three element subsets of groups. \\
       Aequat. Mathem. 22, 140-152 (1981).\\

\item  G.A. Freiman\\
       What is the structure of $K$ if $K + K$ is small?\\
       Lecture Notes in Mathematics 1240, 109-134 (1987).\\

\item  G.A. Freiman, B.M. Schein\\
       Group and semigroup theoretic considerations inspired by
       inverse problems of additive number theory.\\
       Lecture Notes in Mathematics 1320, 121-140 (1988).\\

\item  G.A. Freiman\\
       On extremal additive problems of Paul Erdos.\\
       Ars Combinatoria 26B, 93-114 (1988).\\

\item  G.A. Freiman, P. Erdos\\
       On two additive problems.\\
       Journal of Number Theory 24, 1-12 (1990).\\

\item  G.A. Freiman\\
       Subset-sum problem with different summands.\\
       Congressus Numerantium 70, 207-215 (1990).\\

\item  G.A. Freiman, B.M. Schein\\
       Interconnections between the structure theory of set addition
       and rewritability in groups.\\
       Proc. of Amer. Math. Soc. 113, n.4, 899-910 (1991).\\

\item  G.A. Freiman, B.M. Schein\\
       Structure of $R(3,3)$-groups.\\
       Israel Journal of Mathematics 77, 17-31 (1992).\\

\item  G.A. Freiman\\
       New analytical results in subset sum problem.\\
       Discrete Mathematics 114, 205-218 (1993).\\

\item G. Freiman\\
      On solvability of a system of two boolean linear equations.\\
      Number Theory: New York Seminar, 1991-1995, Springer, 135-150,(1996).\\

\item G. Freiman\\
      Structure theory of set addition.\\
      Structure theory of set addition, Asterisque, 258, 1-33 (1999).\\

\item G.A. Freiman\\
      Structure theory of set addition. II. Results and problems.\\
      Paul Erd\" os and his Mathematics. I, Budapest (Hungary), 1998.\\
      Budapest, 2001, pp. 1-18. Bolyai Society Mathematical Studies X.\\

\item Gregory A. Freiman\\
      On the detailed structure of sets with small additive property.\\
      Combinatorial Number Theory and Additive Group Theory,\\
      Series: Advanced Courses in Mathematics - CRM Barcelona, 233-239 (2009).\\

\item Gregory A. Freiman\\
      Inverse Additive Number Theory XI.\\ 
      Long arithmetic progressions in sets with small sumsets.\\
      Acta Arithmetica, Vol.137, n.4, 325-331, (2009).\\

 \item Ben Green\\
       Approximate groups and their applications: work of Bourgain,Gamburd,
       Helfgott and Sarnak.\\
       arXiv:0911.3354v2 [math.NT],18 Nov 2009.\\
       
\item Ehud Hrushovski\\
      Stable group theory and approximate subgroups.\\
       arXiv.0909.2190v2 [math.LO 18] Sep 2009.\\

\item Renling Jin\\
Freiman's inverse problem with small doubling
property.\\
Adv.\ in Math.,Vol.216,No.2,(2007), 711-725.\\

\item P.Longobardy and M.Maj\\
      The classification of groups with the small squaring property on 3-sets.\\
      Bull.Austral.Math.Soc.,46(1992)263-269\\

\item Martello,S. and P.Toth,
      A Mixture of Dynamic Programming and Branch-and-Bound for the Subset-Sum Problem,
      Management Science,30(1984),765-771\\
          
\item Imre Ruzsa \\
      Arithmetic progressions and the number of sums.\\
      Periodica Math.Hungar.,25:105-111,1992.\\
 
\item Imre Ruzsa \\
       Generalized arithmetic progressions and sum sets.\\
       Acta Math.Hungar.,65(4):379-388,1994\\

\item Tom Sanders\\
      Appendix to: ``Roth's theorem on progressions revisited'',\\
      J. Anal. Math. 104 (2008), 155--192 by J. Bourgain. \\
      J. Anal. Math. 104 (2008), 193--206. \\

\item Terence Tao\\
      Product set estimates for non-commutative groups.\\
      Combinatorica 28 (2008),no. 5, 547-594.\\

\end{enumerate}
\end{document}